\documentclass[a4paper,12pt,oneside]{amsart}
\usepackage{graphics}
\usepackage{blindtext}
\usepackage{amsmath,amssymb,amsthm}
\usepackage{graphicx}

\theoremstyle{definition}

\theoremstyle{definition}

\begin{document}
\title{A note on Quartic Diophantine equation $X^4 - Y^4 = R^2 - S^2$}
\author[S.Muthuvel]{S.Muthuvel}
\date{}
\address{Research Scholar, Department of Mathematics, Faculty of Engineering and Technology, SRM Institute of Science and Technology, Vadapalani Campus, No.1 Jawaharlal Nehru Salai, Vadapalani, Chennai-600026, Tamil Nadu, India}
\email{muthushan15@gmail.com, ms3081@srmist.edu.in}

\author[R.Venkatraman]{R.Venkatraman}
\date{}
\address{Assistant Professor, Department of Mathematics, Faculty of Engineering and Technology, SRM Institute of Science and Technology, Vadapalani Campus, No.1 Jawaharlal Nehru Salai, Vadapalani, Chennai-600026, Tamil Nadu, India}
\email{venkatrr1@srmist.edu.in}  
  \maketitle

\let\thefootnote\relax
\footnotetext{2020 \textit{Mathematics Subject Classification.} Primary: 11D25, 11D45.} %%%%%%%%%%
\footnotetext{\textit{Key words and phrases.} Quartic Diophantine equation; Integer solution; Elementary method.}

\begin{abstract}
 In this paper, we deal with the quartic Diophantine equation $X^4 - Y^4 = R^2 - S^2$ to present its infinitely many integer solutions.
\end{abstract}

\section{Introduction}\label{Sec:1}
In the general formulation of the Diophantine equation given by $x^{n}+y^{n}=u^{n}+v^{n}$, where $n \in \mathbb{Z}$, the specific case when $n=2$ has been elucidated in \cite{Davies,Kersey book,Pasternak}. For $n=4$, the parametric solutions to the aforementioned equation are described in \cite{Choudhry-1, Euler ppr, Hardy book}. More general Diophantine equations with more variables or with integer coefficients that are not all equal to one were taken into consideration by several authors \cite{Choudhry-h, Elkies, Izadi & Nabardi, Janfada & Nabardi}. The authors provided infinite positive integral solutions for various powers, as illustrated in \cite{Babic & Nabardi}.
The objective of this study is to obtain infinitely many integral solutions of
\begin{eqnarray}
X^4 - Y^4=R^2 - S^2 \label{1}
\end{eqnarray}
using each of the specified parametric methods.
In \cite{Baghlaghdam & Izadi form}, the equation
\begin{eqnarray}
a\left(X_{1}^{'5}+X_{2}^{'5}\right)+\sum_{i=0}^{m}a_{i}X_{i}^{5}=b\left(Y_{1}^{'3}+Y_{2}^{'3}\right)+\sum_{i=0}^{n} b_{i}Y_{i}^{3} \label{s1}
\end{eqnarray} 
is examined, where $m,n \in \mathbb{N} \cup \{0\}, a,b \neq 0 \, \text{and} \, a_{i},b_{i}$ are fixed arbitrary rational numbers.  The solution to \eqref{s1} is found by transforming the equation into either a cubic or a quartic elliptic curve with a positive rank, and is determined through the application of elliptic curve theory.
The authors demonstrated that in \big[\cite{Baghlaghdam & Izadi high}, Main Theorem 2\big], the equation

$$\sum_{i=1}^{n}p_{i} x_{i}^{a_i}= \sum_{j=1}^{m}q_{j} y_{j}^{b_j}$$

\noindent where $m,n,a_{i},b_{j}\in\mathbb{N}, \, p_{i},q_{j} \in \mathbb{Z}, \, i=1,2,\dots,n, \, j=1,2,\dots,m$ 
 has a parametric solution and infinitely many solutions in non-zero integers under the condition that either there exists an $i$ such that $p_i = 1$
and $(a_{i},a_{1}\dots a_{i-1}a_{i+1}\dots a_{n}b_{1}\dots b_{m}) = 1$ or there exists a $j$ such that $q_j = 1$ and
$(b_{j},a_{1}\dots a_{n}b_{1}\dots b_{m}) = 1$. While linear transformations are utilized in this article, we put forth an alternative strategy along with distinct conditions for the integer coefficients to solve Eq. \eqref{1}.

\section{Solving the Diophantine equation $X^4 - Y^4 = R^2 - S^2$}
The trivial solution of the equation \eqref{1} is $(X,Y,R,S)=(m,n,m^2 , n^2)$ for $m, n \in \mathbb{Z}$. We consider three different linear transformations, each of which yields a distinct set of infinitely many integer solutions for Eq.\eqref{1}.
\subsection{Method-1}
Consider the linear transformations,
\begin{eqnarray}
  X=px+u, \qquad Y=qx-u, \qquad R=x+v, \qquad S=px+v \label{M1LT}
\end{eqnarray}
where $p,q,u,v \in \mathbb{Z}$. Introducing \eqref{M1LT} in \eqref{1}, we get
\begin{eqnarray}
   \alpha x^4 + \beta x^3 + \gamma x^2 + \delta x &=& 0 \label{M1_1}
\end{eqnarray}
\noindent where
\begin{align}
  \begin{aligned}
   \alpha &= p^4 - q^4, \\       \gamma &= 6p^2 u^2 - 6q^2 u^2 + p^2 - 1,
  \end{aligned}
  &&
  \begin{aligned}
   \beta &= 4p^3 u + 4q^3 u, \\       \delta &= 4pu^3 + 4qu^3 - 2v + 2pv \label{M1_2}
  \end{aligned}
 \end{align}
\noindent For $\delta=0$ in \eqref{M1_2}, we obtain
$$(2p+2q)u^3 = v(1-p)$$
Further, we put $u=t, v=t^3$ and get $p=\frac{1-2q}{3}$. \\ In \eqref{M1_2}, equating the like terms $\gamma = 0$
   $$\left(q + 1\right)\left[\left(-15t^2 +2\right)q + \left(3t^2 -4\right)\right]= 0$$
\noindent On simplifying the above expression, we have $q=\frac{4 - 3t^2}{2 - 15t^2}$ and therefore \eqref{M1_1} becomes
\begin{eqnarray*}
  \alpha x^4 + \beta x^3 &=& 0 \\
   x&=& \frac{t\left(-405t^8 + 459t^6 - 1404t^4 + 600t^2 - 56\right)}{3\left(27t^6 - 27t^4 + 36t^2 - 10\right)}
  \end{eqnarray*}
\noindent Plugging the values $p,q,u,v$ in \eqref{M1LT}, we acquire
\begin{eqnarray}
\begin{aligned}
  X &= \frac{t\left(1215t^{10} - 1782t^8 + 4671t^6 - 774t^4 - 366t^2 + 52\right)}{3(2-15t^2)\left(27t^6 - 27t^4 + 36t^2 - 10\right)} \\
  Y &= \frac{t\left(1215t^{10} - 1782t^8 + 4671t^6 - 5634t^4 + 1902t^2 - 164\right)}{3(2-15t^2)\left(27t^6 - 27t^4 + 36t^2 - 10\right)} \\
  R &= \frac{-t\left(324t^8 - 378t^6 + 1296t^4 - 570t^2 + 56\right)}{3\left(27t^6 - 27t^4 + 36t^2 - 10\right)} \\
  S &= \frac{t\left(810t^8 + 1512t^6 + 1674t^4 - 1092t^2 + 112\right)}{3(2-15t^2)\left(27t^6 - 27t^4 + 36t^2 - 10\right)}
\end{aligned}
\end{eqnarray}
\noindent After eliminating the denominators from the above equations, we have
\begin{eqnarray*}
\begin{aligned}
  X &= 1215t^{11} - 1782t^9 + 4671t^7 - 774t^5 - 366t^3 + 52t \\
  Y &= 1215t^{11} - 1782t^9 + 4671t^7 - 5634t^5 + 1902t^3 - 164t \\
  %R &= 3t(2-15t^2)^2 \left(27t^6 - 27t^4 + 36t^2 - 10\right)\left(324t^8 - 378t^6 + 1296t^4 - 570t^2 + 56\right) \\
R&=5904900t^{19}-14368590t^{17}+41898546t^{15}-55842858t^{13}+58236894t^{11} \\ 
& \qquad -36547200t^9+12314916t^7-2186784t^5+193392t^3-6720t \\
  %S &= 3t(15t^2-2)\left(27t^6 - 27t^4 + 36t^2 - 10\right)\left(810t^8 + 1512t^6 + 1674t^4 - 1092t^2 + 112\right)
S&=984150t^{17}+721710t^{15}+1395306t^{13}-1476954t^{11}+3664440t^9-3124332t^7 \\
& \qquad +1027296t^5-140112t^3+6720t
\end{aligned}
\end{eqnarray*}
Thus, we get an integer solution $(X,Y,R,S)$ for equation \eqref{1} for every $t \in \mathbb{Z}$. Hence,
the presented method generates infinitely many integer solutions for the initial
equation \eqref{1}.
\subsection{Method-2}
Introducing a new transformation in \eqref{1},
\begin{eqnarray}
  X=px+u, \qquad  Y=qx-u, \qquad  R=x+v, \qquad  S=px-v \label{M2LT}
\end{eqnarray}
\noindent $p,q,u,v \in \mathbb{Z}$, leads us to the equation of the form
\begin{eqnarray}
  Ax^4 + Bx^3 + Cx^2 + Dx &=& 0 \label{M2_1}
\end{eqnarray}
\noindent where
\begin{align}
  \begin{aligned}
   A &= p^4-q^4, \\       C &= 6p^2 u^2 - 6q^2 u^2 + p^2 - 1,
  \end{aligned}
  &&
  \begin{aligned}
   B &= 4p^3 u + 4q^3 u, \\       D &= 4pu^3 + 4qu^3 - 2v - 2pv \label{M2_2}
  \end{aligned}
 \end{align}
\noindent For $D=0$ in \eqref{M2_2}, we get
$$(2p+2q)u^3 = v(1+p)$$
Further, we set $u=t, v=t^3$ and attain $p=1-2q$. In \eqref{M2_2}, equating the like terms $C = 0$
$$\left(q - 1\right)\left[\left(9t^2 +2\right)q - 3t^2\right]= 0$$
\noindent Using the above equation, we obtain $q=\frac{3t^2}{9t^2 + 2}$ and therefore \eqref{M2_1} becomes
\begin{eqnarray*}
  A x^4 + B x^3 &=& 0 \\
   x&=& \frac{-t\left(243t^8 + 297t^6 + 216t^4 + 72t^2 +8\right)}{\left(27t^6 + 27t^4 + 12t^2 + 2\right)}
  \end{eqnarray*}
\noindent Applying the values $p,q,u,v$ in \eqref{M2LT}, we acquire
\begin{eqnarray}
\begin{aligned}
  X &= \frac{-3t\left(27t^8 + 36t^6 + 27t^4 + 12t^2 + 2\right)}{\left(27t^6 + 27t^4 + 12t^2 + 2\right)} \\
  Y &= \frac{-t\left(81t^8 + 108t^6 + 81t^4 + 24t^2 + 2\right)}{\left(27t^6 + 27t^4 + 12t^2 + 2\right)} \\
  R &= \frac{-2t\left(108t^8 + 135t^6 + 102t^4 + 35t^2 + 4\right)}{\left(27t^6 + 27t^4 + 12t^2 + 2\right)} \\ \label{M2_F}
  S &= \frac{-2t\left(54t^8 + 81t^6 + 60t^4 + 25t^2 + 4\right)}{\left(27t^6 + 27t^4 + 12t^2 + 2\right)}
\end{aligned}
\end{eqnarray}
\noindent By cancelling the denominators in \eqref{M2_F},
\begin{eqnarray*}
\begin{aligned}
  %X &= 3t\left(27t^6 + 27t^4 + 12t^2 + 2\right)\left(27t^8 + 36t^6 + 27t^4 + 12t^2 + 2\right) \\
X&=2187t^{15}+5103t^{13}+6075t^{11}+4617t^9+2322t^7+756t^5+144t^3+12t \\
  %Y &=  t\left(27t^6 + 27t^4 + 12t^2 + 2\right)\left(81t^8 + 108t^6 + 81t^4 + 24t^2 + 2\right)\\
Y&=2187t^{15}+5103t^{13}+6075t^{11}+4293t^9+1890t^7+504t^5+72t^3+4t \\
  %R &= 2t \left(27t^6 + 27t^4 + 12t^2 + 2\right)^3 \left(108t^8 + 135t^6 + 102t^4 + 35t^2 + 4\right)\\
R&=4251528t^{27}+18068994t^{25}+38381850t^{23}+52986636t^{21}+52435512t^{19}\\
& \qquad+38959218t^{17}+22221378t^{15}+9803592t^{13}+3331692t^{11}+858168t^9 \\
& \qquad +162216t^7+21216t^5+1712t^3+64t \\
  %S &= 2t \left(27t^6 + 27t^4 + 12t^2 + 2\right)^3 \left(54t^8 + 81t^6 + 60t^4 + 25t^2 + 4\right)
S&=2125764t^{27}+9565938t^{25}+21139542t^{23}+30154356t^{21}+30784212t^{19}\\
& \qquad +23632722t^{17}+13977846t^{15}+6426864t^{13}+2290356t^{11}+623160t^9 \\
& \qquad +125496t^7+17664t^5+1552t^3+64t
\end{aligned}
\end{eqnarray*}
For any $t \in \mathbb{Z}$, we obtain an integer solution $(X,Y,R,S)$ to equation \eqref{1}. Consequently, the proposed method yields an infinite number of integer solutions to the starting equation \eqref{1}.
\subsection{Method-3}
In this method, we deal with different transformations in \eqref{1}. Let
\begin{eqnarray}
  X=v, \qquad  Y=px+v, \qquad  R=qx+u, \qquad  S=x+u \label{M3LT}
\end{eqnarray}
\noindent $p,q,u,v \in \mathbb{Z}$. In subsection-1, by introducing these linear transformations in \eqref{1}, leads us to the equation of the form
\begin{eqnarray}
  ax^4 + bx^3 + cx^2 + dx &=& 0 \label{M3_1}
\end{eqnarray}
\noindent where
\begin{align}
  \begin{aligned}
   a &= p^4, \\       c &= 6p^2 v^2 + q^2 - 1,
  \end{aligned}
  &&
  \begin{aligned}
   b &= 4p^3 v, \\       d &= 4pv^3 - 2u + 2qu  \label{M3_2}
  \end{aligned}
 \end{align}
\noindent For $d=0$ in \eqref{M3_2}, we obtain
$$(2p)v^3 = u(1 - q)$$
Additionally, we put $u=t^3, v=t$ and get $p=\frac{1 - q}{2}$. In \eqref{M3_2}, equating the like terms $c = 0$
$$\left(q - 1\right)\left[\left(3t^2 +2\right)q - \left(3t^2 - 2\right)\right]= 0$$
\noindent Thus, we get $q=\frac{3t^2-2}{3t^2 + 2}$ and therefore \eqref{M3_1} becomes
\begin{eqnarray*}
  a x^4 + b x^3 &=& 0 \\
    x&=& \frac{-2t\left(81t^8 + 216t^6 + 216t^4 + 96t^2 + 16\right)}{\left(27t^6 + 54t^4 + 36t^2 + 8\right)}
  \end{eqnarray*}
\noindent Taking the values $p,q,u,v$ and applying it in \eqref{M3LT}, we get
\begin{eqnarray}
\begin{aligned}
  X &= t \\
  Y &= \frac{-3t\left(81t^8 + 216t^6 + 216t^4 + 96t^2 + 16\right)}{(3t^2 + 2)\left(27t^6 + 54t^4 + 36t^2 + 8\right)} \\
  R &= \frac{-t\left(405t^{10} + 756t^8 + 216t^6 - 384t^4 - 304t^2 - 64\right)}{(3t^2 + 2)\left(27t^6 + 54t^4 + 36t^2 + 8\right)} \\
  S &= \frac{-t\left(135t^8 + 358t^6 + 396t^4 + 184t^2 + 32\right)}{\left(27t^6 + 54t^4 + 36t^2 + 8\right)}
\end{aligned}
\end{eqnarray}

\noindent Neglecting the denominators from the above equation,
\begin{eqnarray*}
\begin{aligned}
  %X &= t^2(3t^2 + 2)\left(27t^6 + 54t^4 + 36t^2 + 8\right) \\
X&= 81t^{10}+216t^8+216t^6+96t^4+16t^2 \\
  %Y &= 3t^2 \left(81t^8 + 216t^6 + 216t^4 + 96t^2 + 16\right)\\
Y&=243t^{10}+648t^8+648t^6+288t^4+48t^2 \\
  %R &= t^3(3t^2 + 2)\left(27t^6 + 54t^4 + 36t^2 + 8\right) \left(405t^{10} + 756t^8 + 216t^6 - 384t^4 - 304t^2 - 64\right) \\
R&=32805t^{21}+148716t^{19}+268272t^{17}+217728t^{15}+18144t^{13}-120960t^{11} \\
& \qquad -112896t^9-49152t^7-11008t^5-1024t^3 \\
  %S &= t^3(3t^2 + 2)\left(27t^6 + 54t^4 + 36t^2 + 8\right) \left(135t^8 + 358t^6 + 396t^4 + 184t^2 + 32\right)
 S&=32805t^{21}+196344t^{19}+532008t^{17}+849312t^{15}+874656t^{13}+600000t^{11} \\ 
&\qquad +273536t^9+79872t^7+13568t^5+1024t^3
\end{aligned}
\end{eqnarray*}
For each $t \in \mathbb{Z}$, equation \eqref{1} has an integer solution $(X,Y,R,S)$. The resulting method produces an infinite number of integer solutions to the initial equation \eqref{1}.

\end{document}